\documentclass{amsart}
\usepackage{amscd} 
\usepackage{amsfonts} 
\usepackage{amssymb} 
\usepackage{latexsym} 
\input{diagrams}
\newcommand{\ncm}{\newcommand}

\def\C{\mathbb{C}\,}

\newtheorem{theorem}{Theorem}[section]
\newtheorem{prop}[theorem]{Proposition}
\newtheorem{lemma}[theorem]{Lemma}
\newtheorem{cor}[theorem]{Corollary}
\newtheorem{prob}[theorem]{Problem}
\newtheorem{lem&def}[theorem]{Lemma \& Definition}
\newtheorem{definition}[theorem]{Definition}

\ncm{\End}{\mbox{\rm End}\,}

\def\Hom{\mbox{\rm Hom}\,}

\def\Im{\mbox{\rm Im}\,}

\def\id{\mbox{\rm id}}

\def\into{\hookrightarrow}
\def\to{\rightarrow}

\def\o{\otimes}    %tensor product 
\def\bra{\langle}
\def\ket{\rangle}

\ncm{\rarr}[1]{\stackrel{#1}{\longrightarrow}}
\ncm{\larr}[1]{\stackrel{#1}{\longleftarrow}}
\def\cop{\Delta}

\def\eps{\varepsilon}

\def\du1{\hat 1}

\def\-1{_{(-1)}}
\def\0{_{(0)}}
\def\1{_{(1)}}
\def\2{_{(2)}}
\def\3{_{(3)}}

\def\|{\, | \,}

\def\du1{\hat 1}

\begin{document}

\title[Galois theory for bialgebroids and depth two]{Galois theory for bialgebroids, depth two and normal Hopf subalgebras}
\author{Lars Kadison}
\address{Matematiska Institutionen \\ G{\" o}teborg 
University \\ 
S-412 96 G{\" o}teborg, Sweden} 
\email{lkadison@c2i.net} 
\date{}
\thanks{My thanks to Tomasz Brzezi\'nski, Alexandre Stolin and Korn\'el Szlach\'anyi 
for discussions related to this paper,
and to Marit for TLC.}
\date{} 

\begin{abstract} 
We reduce certain proofs in \cite{KS, Karl, LK} to depth two quasibases from one side only, a  minimalistic approach  which leads to a characterization of Galois extensions for finite
projective bialgebroids without the Frobenius extension property.
We prove that 
a proper algebra extension is a left $T$-Galois extension for some right finite projective
left bialgebroid over some algebra $R$ if and only if it is a left depth two and left balanced
extension. Exchanging left and right in this statement,
we have a characterization of right Galois extensions for
left finite projective right bialgebroids.  Looking to examples of depth two, 
we establish that a Hopf subalgebra is normal if and only if
it is a Hopf-Galois extension.  
We characterize finite weak Hopf-Galois extensions using an alternate Galois 
canonical mapping with several corollaries: that these are depth two and
that surjectivity of the Galois mapping implies its bijectivity.    
\end{abstract} 
\maketitle

\section{Introduction and Preliminaries}

Hopf algebroids arise as the endomorphisms of fiber functors from certain tensor categories
to a bimodule category over a base algebra. For example, 
Hopf algebroids over a one-dimensional base algebra are Hopf algebras
while Hopf algebroids over a separable $K$-algebra base are  weak Hopf algebras.  
 Galois theory for right or left bialgebroids were  recently introduced in \cite{Karl, LK, nov} based on  
the theory of Galois corings \cite{BW} and ordinary  definitions of Galois extensions \cite{Mo, CDG} 
with applications to depth two extensions.  In particular,  Frobenius extensions that are right Galois over
a left finite projective right bialgebroid are characterized in \cite{LK} as being of depth two and right balanced.  
Then a Galois theory for Hopf algebroids, especially of Frobenius type, was introduced
in \cite{Bo, BaS} with applications to Frobenius extensions of depth two and weak Hopf-Galois extensions over
finite dimensional quantum groupoids.
Although they break with the tradition of defining Galois extensions over bialgebras
and have a more complex definition, Galois extensions over Hopf algebroids have more properties in common 
with Hopf-Galois extensions.  However, several of these properties will follow from any Galois theory for bialgebroids
which is in possession of two Galois mappings equivalent due to a bijective antipode, sometimes denoted by $\beta$ and $\beta'$, as is the case for finite Hopf-Galois theory \cite[ch.\ 8]{Mo}, finite weak Hopf-Galois theory (see the last section in this paper),
 possibly some future, useful weakening of Hopf-Galois theory to Hopf algebroids over a symmetric algebra, a Frobenius
algebra or 
some other type of base algebra.  

In \cite{BaS} a characterization similar to that in \cite{LK} for depth two Frobenius one-sided balanced extensions is 
given in terms of Galois extensions over Hopf algebroids with integrals. This shows in a way that
the main theorem in \cite{LK} makes no essential use of the hypothesis of Frobenius extension
(only that a Frobenius extension is of left depth two if and only if it is of right depth two),
and it would be desirable to remove the Frobenius hypothesis.  This is then the objective of section~2 of this
paper:  to show that Galois extensions over one-sided finite projective bialgebroids are characterized by
one-sided depth two and balance conditions on the extension (Theorem~\ref{th-characterization}). 
This requires among other things some care in re-doing the two-sided arguments in \cite{KS} to show that the structure $ T := (A \o_B A)^B$
on a one-sided depth two extension $A \| B$ with centralizer $R$ is still a one-sided finite projective
right bialgebroid (proposition~\ref{prop-tee}). This provides the objective
of the rest of this section; in the appendix in section~5, we include some related
results for the $R$-dual bialgebroid $S:= \End {}_BA_B$ of
a one-sided depth two extension $A \| B$.
These two sections may be read as an introduction to depth two
theory.  

 A depth two complex subalgebra is a generalization of normal subgroup \cite{KK}.
The question was then raised whether depth two Hopf subalgebras
are precisely the normal Hopf subalgebras ($\supseteq$ in \cite{KK}).  
In a very special case, this is true when the notion of
depth two is narrowed to H-separability \cite{Karl}, an exercise in going up and down with ideals as in commutative
algebra. We  study in section~3 the special case of depth two represented by finite Hopf-Galois extensions:
we show that a Hopf-Galois Hopf subalgebra is normal using a certain algebra epimorphism from the Hopf overalgebra
to the Hopf algebra which is coacting Galois, and comparing dimensions of the kernel with the associated Schneider coalgebras. 

A special case of Galois theory for bialgebroids is weak Hopf-Galois theory \cite{BW,CDG,Karl, LK},
(where Hopf-Galois theory is in turn a special case):  for depth two extensions, each type of Galois extension
occurs as we move from any centralizer to separable centralizers to one-dimensional centralizers.  
Conversely, each type of Galois extension, so long as it is finitely generated, is of depth two \cite{KS, Karl, LK}.
In section~4 we complete the proof that a weak Hopf-Galois extension is left  depth two
by studying the alternative Galois mapping $\beta': A \o_B A \to A\o H$
where $\beta'(a \o a') = a\0 a' \o a\1$. As a corollary we find an interesting factorization of the Galois isomorphism
of a weak Hopf algebra over its target subalgebra.  In a second corollary,  a direct proof is given that a surjective
Galois mapping for an  $H$-extension
 is automatically bijective, if $H$ is a finite dimensional weak Hopf algebra.   Finally, it is shown by somewhat
different means than in \cite{BW} that a weak bialgebra in Galois extension of
its target subalgebra has an antipode reconstructible from the Galois mapping.  We provide some evidence
for more generally a
weak bialgebra, which coacts Galois on an algebra over a field, 
 having an antipode, something which is true for bialgebras by a result of
Schauenburg \cite{Sch}.

Let $K$ be any commutative ground ring in this paper. All algebras
are unital associative $K$-algebras and modules over these are symmetric unital $K$-modules. 
We say that $A \| B$ is an \textit{extension} (of algebras) if there is an algebra homomorphism $B \to A$, \textit{proper} if this
is monic.  This homomorphism induces the natural bimodule
structure ${}_BA_B$ which is most important to our set-up.
The extension $A \| B$ is \textit{left depth two} (left D2) if 
the tensor-square $A \o_B A$ is centrally projective w.r.t.\
$A$ as natural $B$-$A$-bimodules: i.e., 
$$ {}_BA \o_B A_A \oplus * \cong \oplus^n {}_BA_A. $$
This last statement postulates the existence then of a
split $B$-$A$-epimorphism from a direct sum 
of $A$ with itself $n$ times to $A \o_B A$.  

Making the clear-cut identifications  $\Hom ({}_BA \o_B A_A, {}_BA_A)
\cong $ $\End {}_BA_B$ and $\Hom ({}_BA_A, {}_BA \o_B A_A) \cong $
$(A \o_B A)^B$, we see that left D2 is characterized by
there being a left D2 quasibase $t_i \in (A \o_B A)^B$
and $\beta_i \in \End {}_BA_B$ such that for all $a, a' \in A$
\begin{equation}
\label{eq: left D2 quasibase}
a \o_B a' = \sum_{i=1}^n t_i \beta_i(a)a'.
\end{equation}
The algebras $\End {}_BA_B$ and $(A \o_B A)^B$ (note
that the latter is isomorphic to $\End {}_A A \o_B A_A$
and thus receives an algebra structure) are so important
in depth two theory that we fix (though not unbendingly) brief notations for these:
$$ S := \End {}_BA_B \ \ \ \ \ \ \  T := (A \o_B A)^B. $$

Similarly, a right depth two extension $A \| B$ is defined by switching from the natural $B$-$A$-bimodules in the definition
above to the natural
$A$-$B$-bimodules on the same structures.  Thus an extension
$A \| B$ is \textit{right D2} if ${}_A A \o_B A_B \oplus * \cong \oplus^m {}_AA_B$.  Equivalently, if there are
$m$ paired elements $u_j \in T$, $\gamma_j \in S$ such that
\begin{equation}
\label{eq: right D2 quasibase}
a \o a' = \sum_{j=1}^m a \gamma_j(a') u_j
\end{equation}
for all $a, a' \in A$.  

A depth two extension is one that is both left and right D2.
These have been studied in \cite{KS,Karl,LK} among others, 
but without a focus on left or right D2 extensions.
Note that
a left D2 extension $A \| B$ has right D2 extension
$A^{\rm op} \| B^{\rm op}$
when we pass to opposite algebras. This gives in
fact a natural one-to-one correspondence between left D2 extensions
and right D2 extensions.

Let $t, t'$ be elements in $T$, where we write $t$
in terms of its components using a notation that
suppresses a possible summation in $A \o_B A$:
$t = t^1 \o t^2$. 
Then the algebra structure on $T$ is simply
\begin{equation}
tt' = {t'}^1 t^1 \o t^2 {t'}^2, \ \ \ \ \ 1_T = 1_A \o 1_A
\end{equation}

There is a standard ``groupoid'' way to produce right and left
bialgebroids, which we proceed to do for $T$.  
There are two commuting embeddings of $R$ and its opposite algebra
in $T$.  A ``source'' mapping $s_R: R \to T$ given by
$s_R(r) = 1_A \o r$, which is an algebra homomorphism.
And a ``target'' mapping $t_R: R \to T$ given by
$t_R(r) = r \o 1_A$ which is an algebra anti-homomorpism
and clearly commutes with the image of $s_R$.  Thus it makes
sense to give $T$ an $R$-$R$-bimodule structure via $s_R$, $t_R$
from the right:  $r \cdot t \cdot r' = $ $t s_R(r') t_R(r) = $
$t(r \o r') = $ $rt^1 \o t^2 r'$, i.e., ${}_RT_R$ is given by
\begin{equation}
\label{eq: bimod}
r\cdot t^1 \o t^2 \cdot r' = rt^1 \o t^2 r'
\end{equation}

\begin{prop}
\label{prop-tee}
Suppose $A \| B$ is either a right or a left D2 extension.
Then $T$ is a right $R$-bialgebroid, which is either
left f.g.\ $R$-projective or right f.g.\ $R$-projective
respectively.
\end{prop}
\begin{proof}
First we suppose $A \| B$ is left D2 with quasibases
$t_i \in T$, $\beta_i \in S$.  The proof that $T$ is
a right $R$-bialgebroid in \cite[5.1]{KS} carries through
verbatim except in one place where a right D2 quasibase made
a brief appearance, where coassociativity of the coproduct
needs to be established through the introduction of an isomorphism.  Thus we need to see that
$$ T \o_R T \o_R T \stackrel{\cong}{\longrightarrow} 
(A \o_B A \o_B A \o_B A)^B$$
via $t \o t' \o t'' \mapsto t^1 \o t^2 {t'}^1 \o {t'}^2 {t''}^1 \o {t''}^2$.  
The  inverse is given by
$$ a_1 \o a_2 \o a_3 \o a_4 \longmapsto 
\sum_{i,j} t_i \o_R t_j \o_R (\beta_j(\beta_i(a_1)a_2)a_3 \o_B a_4). $$
for all $a_i \in A$ ($i = 1,2,3,4$).

In the case that we only use a right D2 quasibase,
this inverse is given by
\begin{equation}
\label{eq: coassoc}
a_1 \o a_2 \o a_3 \o a_4 \longmapsto
\sum_{j,k} a_1 \o a_2 \gamma_k(a_3 \gamma_j(a_4)) \o_R u_k \o_R u_j.
\end{equation}
Both claimed inverses are easily verified as such by using the
right and left D2 quasibase equations repeatedly.

The module $T_R$ is finite projective since eq.~(\ref{eq: left D2
quasibase}) implies a dual bases equation $t =$ $ \sum_i t_i f_i(t)$, 
for each $t \in T \subseteq A\o_B A$, 
 where $f_i(t) := \beta_i(t^1)t^2$ define $n$ maps in
$\Hom (T_R, R_R)$.

Suppose $A \| B$ is right D2 with quasibase $u_j \in T$, $\gamma_j \in
S$.  The algebra structure on $T$ is given in the introduction
above as is the $R$-$R$-bimodule structure.
What remains is specifying the $R$-coring structure on $T$
and checking the five axioms of a right bialgebroid.
The coproduct $\cop: T \to T \o_R T$ is given
by 
\begin{equation}
\label{eq: cop}
\cop(t) := \sum_j (t^1 \o_B \gamma_j(t^2)) \o_R u_j,
\end{equation}
 which is clearly left $R$-linear, and right $R$-linear as well 
since
$$ \cop(tr) = \sum_j t^1 \o \gamma_j(t^2 r) \o u_j \stackrel{\cong}{\longmapsto} t^1 \o 1 \o t^2 r $$
under the isomorphism $T \o_R T \cong (A \o_B A \o_B A)^B$
given by $t \o t' \mapsto $ $t^1 \o t^2 {t'}^1 \o {t'}^2$,
which is identical to the image of 
$$ \cop(t)r = \sum_j t^1 \o \gamma_j(t^2) \o u_j r \mapsto 
t^1 \o 1 \o t^2 r .$$

Coassociativity $(\cop \o \id_T) \cop = (\id_T \o \cop) \cop$ follows
from applying the isomorphism $$T \o_R T \o_R T \cong
(A \o_B A \o_B A \o_B A)^B$$ given above in this proof 
to the left-hand and right-hand sides applied to a $t \in T$: 
$$ \sum_j \cop(t^1 \o \gamma_j(t^2)) \o u_j = \sum_{j,k}
t^1 \o \gamma_k(\gamma_j(t^2)) \o_R u_k \o_R u_j 
\stackrel{\cong}{\longmapsto} t^1 \o_B 1_A \o_B 1_A \o_B t^2. $$
$$ \sum_j (t^1 \o \gamma_j(t^2)) \o_R \cop(u_j)
= \sum_{j,k} (t^1 \o \gamma_j(t^2)) \o_R (u_j^1 \o \gamma_k(u_j^2))
\o_R u_k $$
which also maps into $t^1 \o_B 1_A \o_B 1_A \o_B t^2$ under the
same isomorphism.

The counit $\eps: T \to R$ of the $R$-coring $T$ is given by
\begin{equation}
\eps(t) := t^1 t^2
\end{equation}
i.e., the multiplication mapping $A \o_B A \to A$ restricted
to $T$ (and taking values in $R$ since $bt = tb$ for all $b \in B$).
Clearly, $\eps(rtr') = r \eps(t) r'$ for $r,r' \in R$, $t \in T$,
and that $(\id_T \o_R \eps)\cop = \id_T$ $= (\eps \o_R \id_T)\cop$
since $\sum_j t^1 \gamma_j(t^2)u_j = t$
and $\sum_j \gamma_j(a)u_j^1u_j^2 = a$ for $a \in A$. 
This shows $(T,R,\cop, \eps)$ is an $R$-coring.

We next verify the five axioms of a right bialgebroid \cite[2.1]{KS}.

\begin{enumerate}
\item $\cop(1_T) = 1_T \o 1_T$ since $\gamma_j(1_A) \in R$
and $1_T = \sum_j \gamma_j(1_A)u_j$.

\item $\eps(1_T) = 1_A$ since $1_T = 1_A \o 1_A$.

\item $\eps(tt') = \eps(t_R(\eps(t)) t') = \eps(s_R(\eps(t)) t')$
($t,t' \in T$) since $\eps(tt') =$ $ {t'}^1 t^1 t^2 {t'}^2$, 
$t_R(\eps(t)) =$ $ t^1 t^2 \o_B 1_A$ and $s_R(\eps(t)) =$ $ 1_A \o_B
t^1 t^2 $.

\item $(s_R(r) \o 1_T)\cop(t) = (1_T \o t_R(r))\cop(t)$ for
all $r \in R$, $t \in T$ since the left-hand side is
$$\sum_j (t^1 \o_B r \gamma_j(t^2)) \o_R u_j \stackrel{\cong}{\longmapsto} t^1 \o_B r \o_B t^2 $$
under the isomorphism $T \o_R T \cong (A \o_B A \o_B A)^B$
given by $t \o_R t' \mapsto$ $ t^1 \o_B t^2 {t'}^1 \o_B {t'}^2$
and the right-hand side is equal to
$$ \sum_j (t^1 \o_B \gamma_j(t^2)) \o_R (u^1_j r \o_B u^2_j) 
\stackrel{\cong}{\longmapsto} t^1 \o_B r \o_B t^2 $$ 
with the same image element.  

\item $\cop(tt') = \cop(t)\cop(t')$ for all $t,t' \in T$ in the tensor subalgebra (denoted by $T \times_R T$ with the straightforward
tensor multiplication)
of $T \o_R T$ (which makes sense thanks to the previous axiom). 
This follows from both sides having the image element ${t'}^1 t^1 \o 1_A
\o t^2 {t'}^2 $ under the isomorphism $T \o_R T \cong $
$(A \o_B A \o_B A)^B$, which is clear for the left-hand side of the equation
and for  the right-hand side we note it equals 
$$ \sum_{j,k} ({t'}^1 t^1 \o_B \gamma_j(t^2) \gamma_k({t'}^2) )
\o_R (u_k^1 u_j^1 \o_B u_j^2 u_k^2 ). $$
Now apply $t \o t' \mapsto t^1 \o_B t^2 {t'}^1 \o_B {t'}^2$
and the right D2 quasibase equation twice. 
\end{enumerate}
This completes the proof that $(T, R, s_R, t_R, \cop, \eps)$
is a right bialgebroid. 

Finally ${}_RT$ is finite projective via an application
of the right D2 quasibase eq.~(\ref{eq: right D2 quasibase}).
\end{proof}

A right comodule algebra is an algebra in the tensor category of right $R$-comodules
\cite{BaS}.  In detail, the definition is equivalent to the following.   
  
\begin{definition}
Let $T$ be any right bialgebroid $(T, R, \tilde{s}, \tilde{t},$ $ \cop, $ $ \eps)$ over any base algebra $R$.   
A right $T$-comodule algebra $A$ is an algebra $A$ with algebra homomorphism $R \to A$ (providing the $R$-$R$-bimodule structure on $A$)
 together with a coaction $\delta: A \to A \o_{R} T$, where  values $\delta(a)$ are denoted by the Sweedler
notation $a\0 \o a\1$, 
such that $A$ is a right $T$-comodule over the $R$-coring $T$ \cite[18.1]{BW}, 
$ \delta(1_{A}) = 1_{A} \o 1_{T} $, $ra\0 \o a\1 =$ $ a\0 \o \tilde{t}(r)a\1$
for all $r \in R$,
and $\delta(aa') = \delta(a) \delta(a')$ for all $a,a' \in A$.   The subalgebra of coinvariants
is ${A}^{\rm co \, T} := \{ a \in A | \delta(a) = a \o 1_{T} \}$.
We also call $A$ a right $T$-extension of ${A}^{\rm co \, T}$.   
\end{definition} 

\begin{lemma}
For the right $T$-comodule $A$ introduced just above,
 $R$ and ${A}^{\rm co \, T}$ commute in $A$.
\end{lemma}
\begin{proof}
We note that 
$$\delta (rb) = b \o_R \tilde{s}(r) = \delta (br) $$
for $r \in R$, $b \in {A}^{\rm co \, T}$. But $\delta$
is injective by the counitality of comodules, so
$rb = br$ in $A$ (suppressing the morphism $R \to A$).  
\end{proof}

\begin{definition}
Let $T$ be any  right bialgebroid
over any algebra $R$.  
A  $T$-comodule algebra $A$ is a right $T$-Galois extension of its coinvariants
$B$ if the (Galois) mapping 
$\beta: A \o_{B} A \to A \o_{R} T$ defined
by $\beta(a \o a') = a{a'}\0 \o {a'}\1 $ is bijective.  
\end{definition} 

Left comodule algebras over left bialgebroids and their left Galois extensions
are defined similarly, the details of which are in \cite{nov}.  The values of the coaction is in this case denoted
by $a\-1 \o a\0$ and the Galois mapping by $a \o a' \mapsto a\-1 \o a\0 a'$. When we pass to opposite algebras, we note that a left $T$-Galois extension $A \| B$
has right $T^{\rm op}$-Galois extension $A^{\rm op} \| B^{\rm op}$.

%%%%%%%%%%%%%%%%%%%%%%%%%%%%%%%%%%%%%%%%%%%%%%%%%%%%%%%%%%%%%%%%%%%%
\section{A characterization of Galois extensions for bialgebroids}

We recall that a module ${}_AM$ is \textit{balanced} if all 
the endomorphisms of the natural module $M_E$ where $E = \End {}_AM$
are uniquely left multiplications by elements of $A$:
$A \stackrel{\cong}{\longrightarrow} \End M_E$ via $\lambda$.  
In particular, ${}_AM$ must be faithful.
 
\begin{theorem}
\label{th-characterization}
Let $A \| B$ be a proper algebra extension. Then
\begin{enumerate}
\item  $A \| B$ is a right $T$-Galois extension for
some left finite projective right bialgebroid $T$ over
some algebra $R$ if and only if $A \| B$ is right D2 and right balanced.
\item $A \| B$ is a left $T$-Galois extension for
some right finite projective left bialgebroid $T$ over
some algebra $R$ if and only if $A \| B$ is left D2 and left balanced.
\end{enumerate}
\end{theorem}
\begin{proof}
($\Rightarrow$)
Suppose $T$ is a left finite projective right bialgebroid
over some algebra $R$. 
Since  ${}_RT \oplus * \cong {}_RR^t$ for some positive integer $t$,  we apply to this 
the functor $A \o_R -$ from left $R$-modules into $A$-$B$-bimodules which results in 
${}_AA\! \o_B\! A_B \oplus * \cong {}_AA_B^t$, after using the Galois $A$-$B$-isomorphism 
$A \o_B A \cong  $ $A \o_R T$. Hence, $A | B$ is right D2.

Let $\mathcal{E} := \End A_B$.  We show $A_B$ is balanced by
the following device. 
Let $R$ be an algebra, $M_R$ and ${}_RV$ modules with ${}_RV$ finite projective.  If 
$\sum_j m_j \phi(v_j) = 0$ for all $\phi$ in the left $R$-dual ${}^*V$, then $\sum_j m_j \o_R v_j = 0$.
This follows immediately by using dual bases $f_i \in {}^*V$, $w_i \in V$.  

Given $F \in \End {}_{\mathcal{E}}A$, it suffices to show that $F = \rho_b$ for some
$b \in B$.  Since $\lambda_a \in \mathcal{E}$, $F \circ \lambda_a = \lambda_a \circ F$ for all $a \in A$,
whence $F = \rho_{F(1)}$.  Designate $F(1) = x$.  If we show that $x\0 \o x\1 = x \o 1$ after applying
the right $T$-valued coaction on $A$,
then $x \in A^{\rm co \, T} = B$. For each $\alpha \in \Hom ({}_RT, {}_RR)$, define $\overline{\alpha} \in \End A_B$
by $\overline{\alpha} (y) = y\0 \alpha(y\1)$. Since $\rho_r \in \mathcal{E}$ for each $r \in R$ by lemma, 
$$x \alpha(1_T) = F(\overline{\alpha} (1_A)) = \overline{\alpha} (F(1_A)) = x\0 \alpha(x\1) $$
for all $\alpha \in {}^*T$. Hence $x\0 \o_R x\1 = x \o 1_T$.

($\Leftarrow$) It follows from the proposition
that a right D2 extension $A \| B$ has a left finite projective
right bialgebroid $T := (A \o_B A)^B$ over the centralizer $R$
of the extension. Let $R \into A$ be the inclusion mapping.
We check that $A$ is a right $T$-comodule algebra via the
coaction $\rho_R: A
\to A \o_R T$ on $A$ given by
\begin{equation}
\label{eq: coaction}
\rho_R(a) = a\0 \o a\1 := \sum_j \gamma_j(a) \o u_j.
\end{equation}
First, we demonstrate several properties by using the isomorphism
$\beta^{-1}: A \o_R T \stackrel{\cong}{\longmapsto} A \o_B A$
given by $\beta^{-1}(a \o t) = at = at^1 \o t^2$ \cite[3.12(iii)]{KS}
with inverse $\beta(a \o a')= \sum_j a{a'}\0 \o {a'}\1$
(cf.\ right D2 quasibase eq.~(\ref{eq: right D2 quasibase})). 
This shows straightaway that the Galois mapping $\beta: A \o_B A
\to A \o_R T$ is  bijective.
Then $A \o_R T \o_R T \cong A \o_B A \o_B A$ via $\Phi := (\id_A \o \beta^{-1})(\beta^{-1} \o \id_T)$,
so coassociativity $(\id_A \o \cop)\rho_R = (\rho_R \o \id_T)\rho_R$ follows
from 
$$\Phi( \id \o \cop_T) \circ \rho_R) = \sum_{j,k} \gamma_j(a) u^1_j \o_B \gamma_k(u^2_j)u^1_k \o_B u^2_k =
\sum_k 1 \o \gamma_k(a)u^1_k \o u^2_k = 1 \o 1 \o a
$$
$$= \sum_{j,k} \gamma_k(\gamma_j(a))u^1_k \o u^2_ku^1_j \o u^2_k = 
\Phi( (\rho_R \o \id)\rho_R(a)). $$ 
We note that $\rho_R$ is right $R$-linear, since
$$ \rho_R(ar) = \sum_j \gamma_j(ar) \o u_j \stackrel{\beta^{-1}}{\longmapsto} 1 \o_B ar = \beta^{-1}(\rho_R(a)r) $$
since $\rho_R(a)r = \sum_j \gamma_j(a) \o u_jr$. 
Also, $a\0 \eps_T(a\1) = $ $ \sum_j \gamma_j(a)u^1_j u^2_j  = a$
for all $a \in A$.  

Next, 
$$ \beta^{-1} (r \cdot a\0 \o a\1) = \sum_j r\gamma_j(a)  u_j = r \o_B a =
\sum_j \gamma_j(a)u^1_j r \o u^2_j = \beta^{-1}(a\0 \o t_R(r) a\1). $$
Whence the statement $\rho_R(aa') = \rho_R(a)\rho_R(a')$
makes sense for all $a,a' \in A$. 
We check the statement:  
$$\beta^{-1}(\rho_R(a)\rho_R(a')) = \sum_{j,k} \gamma_j(a) \gamma_k(a') u_j u_k =
 \sum_{j,k} \gamma_j(a) \gamma_k(a') u_k^1 u_j^1 \o u_j^2 u_k^2 
$$
$$ = 1 \o aa' = \sum_j \gamma_j(aa')u_j = \beta^{-1}(\rho_R(aa')). $$
Also $\rho_R(1_A) = 1_A \o_R 1_T$ since $\gamma_j(1_A) \in R$.
Finally we note that for each $b\in B$
$$ \rho_R(b) = \sum_j \gamma_j(b) \o_R u_j = b \o \sum_j \gamma_j(1)u_j = b \o 1_T $$
so $B \subseteq A^{\rm co \, \rho_R}$.  
Conversely, if $\rho_R(x) = x \o 1_T$ $ = \sum_j \gamma_j(x) \o u_j$ applying $\beta^{-1}$ we
obtain $x \o_B 1 = 1 \o_B x$. Let $f \in \End A_B$.   Then
applying $\mu (f \lambda(a) \o \id)$ to this we obtain 
 $f(ax) = f(a)x$  since $\lambda(a) \in \End A_B$
for each $a\in A$. It follows that $f \rho(x) = \rho(x) f$
so $\rho(x) \in \End {}_{\mathcal{E}}A$.  Since $A_B$ is balanced,
$\rho(x) = \rho(b)$ for some $ b\in B$, whence $x = b \in B$.  

The second part of the theorem is proven similarly
(or alternatively, apply the first part with the opposite algebra technique mentioned
in the introduction). In the $\Leftarrow$ direction,
we convert the right $R$-bialgebroid $T$ to a left $R$-bialgebroid $T^{\rm op}$
with $s_L = t_R$, $t_L = s_L$, the same $R$-coring structure and opposite multiplication, which leads to  
the left $R$-$R$-bimodule structure coinciding with the usual $R$-$R$-bimodule structure
on $T$ in eq.~(\ref{eq: bimod}). We then define a left $T^{\rm op}$-comodule algebra structure on $A$
via $\rho_L: A \to T \o_R A$ defined via  left D2 quasibases by  
\begin{equation}
\rho_L(a) = a\-1 \o a\0 := \sum_i t_i \o \beta_i(a).
\end{equation}
The isomorphism $T \o_R A \stackrel{\cong}{\longrightarrow}$ $ A \o_B A$ given
by $t \o a \mapsto t^1 \o t^2a$ is inverse to the Galois mapping 
$\beta_L(a \o a') = $ $a\-1 \o a\0 a'$ by the left D2 quasibase eq.~(\ref{eq: left D2 quasibase}).
One needs the opposite multiplication of $T$ when showing $\rho_L(aa') = \rho_L(a) \rho_L(a')$
for $a,a'\in A$.  
\end{proof}

Let $T$ be a left finite projective right bialgebroid over some algebra $R$ in the next corollary.
 
\begin{cor}
Suppose $A \| B$ is a right $T$-extension.  If the Galois mapping $\beta$ is a split $A$-$B$-monomorphism,
then $A \| B$ is a right $(A \o_B A)^B$-Galois extension.
\end{cor}
\begin{proof}
This follows from ${}_A A \o_B A_B \oplus * \cong {}_A A_B \o_R T$ and the arguments in the first few paragraphs
of the proof above (the balance argument makes only use of $A \| B$ being a right $T$-extension).  
Hence, $A \| B$ is right D2 and right balanced.  Whence $A \| B$ is a right Galois extension w.r.t.
the bialgebroid $(A \o_B A)^B$.
\end{proof}
Notice that $T$ is possibly not isomorphic to $(A \o_B A)^B$.  For example, one might start with 
a Hopf algebra Frobenius extension with split monic Galois map and conclude it is a weak Hopf-Galois extension
(if the centralizer is separable, the antipode being constructible from the Frobenius structure).

%%%%%%%%%%%%%%%%%%%%%%%%%%%%%%%%%%%%%%%%%%%%%%%%%%%%%%%%%%%%%%%%%%%%%%%%%
\section{Galois extended Hopf subalgebras are normal}

There is a question of whether depth two Hopf subalgebras are normal \cite[3.4]{KK}. 
In this section we answer this question in an almost unavoidable
special case, namely, when the Hopf subalgebra forms a Hopf-Galois
extension with respect to the action of a third Hopf algebra.  
Since a depth two extension with one extra condition is a Galois extension for
actions of bialgebroids or weak bialgebras \cite{LK}, 
the situation of
ordinary Hopf-Galois extension would seem to be a critical  step.  

Let $k$ be a field. All Hopf algebras in this section are finite dimensional
algebras over $k$.  
Recall that a Hopf subalgebra $K \subseteq H$ is a Hopf algebra $K$
w.r.t.\  the algebra and coalgebra structure of $H$ (with counit denoted by $\eps$) as well as
stable under the antipode $\tau$ of $H$. Recall the Nichols-Zoeller
result that the natural modules $H_K$ and ${}_KH$ are free.  
$K$ is \textit{normal} in $H$
if $\tau(a\1)x a\2 \in K$ and $a\1 x \tau(a\2) \in K$
for all $x \in K, a \in H$.  Equivalently, if $K^+$ denotes
the kernel of the counit $\eps$, $K$ is a normal Hopf subalgebra
of $H$ if the left algebra ideal and coideal $HK^+$ is equal to the right ideal and coideal $K^+H$
\cite[3.4.4]{Mo}.
   
In considering another special case of D2 Hopf subalgebras,  we showed in \cite{Karl} that H-separable Hopf subalgebras
are normal using favorable properties for H-separable extensions of going down and going up for ideals.  However, we
noted that
such subalgebras are not proper if $H$ is semisimple, e.g.,
$H$ is a complex group algebra.
In \cite[3.1]{KK} we showed that depth two subgroups are normal subgroups
using character theory (for $k = \C$). We also noted the more general
converse that normal Hopf $k$-subalgebras are Hopf-Galois extensions
and therefore D2. Next we extend this to the 
characterization of normal Hopf subalgebras below,  one that we believe 
is not altogether unexpected but unnoted or not adequately exposed in the literature.

  \begin{theorem}
Let $K \subseteq H$ be a Hopf subalgebra.  Then $K$ is normal in
$H$ if and only if $H | K$ is a Hopf-Galois extension.  
\end{theorem}

\begin{proof} ($\Rightarrow$) This is more or less implicit
in \cite[3.4.4]{Mo}, where it is also shown \cite[chs.\ 7,8]{Mo}
 that $H$ is a crossed product by a counital $2$-cocycle of $K$ with 
the quotient Hopf algebra $\overline{H}$ (a cleft $\overline{H}$-extension or
Galois extension with normal basis).
Since $HK^+ = K^+H$ under normality of $K$, it becomes a Hopf ideal, so 
we form the Hopf algebra $\overline{H} := H/HK^+$, which coacts
naturally on $H$ via the comultiplication and quotient projection.
The coinvariants are precisely $K$ since $H_K$ is faithfully flat.
The Galois mapping $\beta: H \o_K H \rightarrow H \o \overline{H}$
given by $\beta(a \o a') = a {a'}\1 \o \overline{{a'}\2}$ is an isomorphism
with inverse given by $x \o \overline{y} \mapsto x \tau(y\1) \o y\2$. 

($\Leftarrow$)  Suppose $H$ is a $W$-Galois extension of $K$ where $W$ is a Hopf algebra
with right coaction $\rho: H \to H \o W$ on $H$.  We define a mapping $\Phi: H \to W$
by $\Phi(h) = \eps_H(h\0) h\1$, i.e., $\Phi = (\eps_H \o \id_W) \circ \rho$.  We note
that $\Phi$ is an algebra homomorphism since $\rho$ and $\eps_H$ are (and augmented since $\eps_W \circ \Phi = \eps_H$). 
Also, $\Phi: H \to W$ is a right $W$-comodule morphism since $H$ is a right $W$-comodule with
$\rho$ and $\cop_W$ obeying a coassociativity rule.  Next we note that $\Phi$ is an epi
since given $w \in W$, there is $\sum_i h_i \o {h'}_i \in H \o_K H$ such that
$1 \o w = \sum_i h_i {{h'}_i}\0 \o {{h'}_i}\1$.  Applying $\eps_H \o \id_W$ to this, we obtain 
$$ \Phi(\sum_i \eps_H(h_i) {h'}_i) = w . $$

We note that $\ker \Phi$ contains $K^+$ since $K = H^{\rm co \, W}$ $ = \{ h \in H \, |\, \rho(h) = h \o 1_W \}$
Consider the coalgebra and right quotient $H$-module $H/ K^+H := \overline{H}$
as well as the coalgebra and left quotient $H$-module $H/ HK^+ := \overline{\overline{H}}$.
In this case, $\Phi$ induces $\overline{\Phi}:
\overline{H} \to W$ and $\overline{\overline{\Phi}}: \overline{\overline{H}} \to W$.
(They are respectively right and left $H$-module morphisms w.r.t.\ the modules $W_{\Phi}$ and ${}_{\Phi}W$.)
By Schneider \cite[1.3]{HJS}, the Galois quotient mapping $\overline{\beta}: H \o_K H \to H \o \overline{H}$
given by $\overline{\beta}(x \o y) = xy\1 \o \overline{y\2}$ is bijective (since $K$ is a left coideal subalgebra of $H$).
But the Hopf subalgebra $K$ is also a right coideal subalgebra satisfying a right-handed version of Schneider's lemma
recorded in \cite[2.4]{FMS}: whence the Galois mapping $\overline{\overline{\beta}}: H \o_K H \to \overline{\overline{H}} \o H$ given by $\overline{\overline{\beta}}(x \o y) = \overline{x\1} \o x\2 y$ is bijective as well.

Observe now that $H_K$ is free of rank $n$, let's say,  so $\beta$ bijective implies that $\dim_k W = n$.
Similarly, $\overline{\beta}$ bijective implies $\dim_k \overline{H} = n$
and $\overline{\overline{\beta}}$ bijective implies $\dim_k \overline{\overline{H}} = n$.  
It follows that the  vector space epimorphisms $\overline{\Phi}: \overline{H} \to W$ and $\overline{\overline{\Phi}}: \overline{\overline{H}} \to W$ are isomorphisms. But $\overline{\Phi}$ factors through $\overline{H} \to H/HK^+H$ induced
by $K^+H \subseteq HK^+H$; similarly, $\overline{\overline{\Phi}}$ factors through $\overline{\overline{H}} \to
H/HK^+H$, so both these canonical mappings are monic.  
It follows that $HK^+ = HK^+H = K^+H$, whence $K$ is a normal Hopf subalgebra in $H$. 
  \end{proof}

In the proof of $\Leftarrow$ above, we can go further to conclude that $\overline{H}$ is a Hopf algebra
isomorphic to $W$ as augmented algebras.  However, the theory of deforming the comultiplication of
a Hopf algebra by a $2$-cocycle
\cite[2.3.4]{Maj} shows that there are pairs of Hopf algebras isomorphic as augmented algebras yet
non-isomorphic as Hopf algebras. Additionally, there are examples of algebra extensions which are Hopf-Galois w.r.t.\ two
different Hopf algebras.  We therefore do not know \textit{a priori}
 if $\overline{H}$ and $W$ are isomorphic as Hopf algebras.

%%%%%%%%%%%%%%%%%%%%%%%%%%%%%%%%%%%%%%%%%%%%%%%%%%%%%%%%%%%%%%%%%%%%%%%%%%%%%%%%%%%%%%%%%%%%%%%%%%%%%%%%%%%%%%%%%%%%%%%
\section{Weak Hopf-Galois extensions are depth two}

In this section we study right Galois extensions of special bialgebroids - the weak Hopf-Galois extensions, cf.\ \cite{BW,CDG, Karl, LK}. By exploiting the antipode in weak Hopf algebras, we find an alternative Galois mapping
which characterizes weak Hopf-Galois extensions.  This leads
to several corollaries that finite weak Hopf-Galois extensions are right as well as left depth two
extensions, that they may be defined by only a surjective Galois map, and that a weak Hopf algebra over
its target separable subalgebra is an example of such.  We propose a number of problems for further study
in the young subject of weak Hopf-Galois extensions.

  Weak Hopf algebras are a special case of Hopf algebroids - those 
with separable base algebra \cite{EN, KS}: the separable
algebra has an index-one Frobenius system which one uses to
convert mappings to the base and tensors over the base
to linear functionals and tensors over a ground field. There is an example of one step in how to conversely
view a weak Hopf algebra $H$ as a Hopf algebroid over its left coideal subalgebra $H^L$ in the proof of corollary~\ref{cor - olary} below.  

Let $k$ be a field. A weak Hopf algebra
$H$ is  first a weak bialgebra, i.e., a $k$-algebra and $k$-coalgebra $(H, \cop, \eps)$
such that the comultiplication $\cop: H \to H \o_k H$ is linear and multiplicative,
$\cop(ab) = \cop(a) \cop(b)$, and the counit is linear
just as for bialgebras;
however, the change (or weakening of the axioms) is that $\cop$ and $\eps$ may not be unital, $\cop(1) \neq
1 \o 1$ and $\eps(1_H) \neq 1_k$, but must satisfy 
\begin{equation}
1\1 \o 1\2 \o 1\3 = (\cop(1) \o 1)(1 \o \cop(1)) =
(1 \o \cop(1))(\cop(1) \o 1)
\end{equation}
 and $\eps$ may not be multiplicative, $\eps(ab) \neq \eps(a) \eps(b)$ but must satisfy $(a,b,c \in H$) 
\begin{equation}
\eps(abc) = \eps(ab\1) \eps(b\2 c) = \eps(a b\2) \eps(b\1 c).
\end{equation}
There are several important projections that result from these axioms:
\begin{eqnarray}
\Pi^L(x) &:=& \eps(1\1 x) 1\2 \\
\Pi^R(x) &:=& 1\1 \eps(x 1\2) \\
\overline{\Pi}^L(x) & := & 1\1 \eps(1\2 x) \\
 \overline{\Pi}^R(x) & := & \eps(x 1\1) 1\2 \ \ \ \ \ (\forall \, x \in H)
\end{eqnarray}
We denote $H^L := \Im \Pi^L = \Im \overline{\Pi}^R$
and $H^R := \Im \Pi^R = \overline{\Pi}^L$. These subalgebras
are separable $k$-algebras \cite{BNS}.

In addition to being a weak bialgebra, a weak Hopf algebra
has an antipode $S: H \to H$ satisfying the axioms
\begin{eqnarray}
S(x\1)x\2 & = & \Pi^R(x) \label{eq: bns2} \\
x\1 S(x\2) & = & \Pi^L(x) \label{eq: bns1} \\
S(x\1) x\2 S(x\3) & = & S(x) \label{eq: bnsa} \ \ \ \  (\forall \, x \in H)
\end{eqnarray}
The antipode turns out to be bijective for finite dimensional
weak Hopf algebras (which we will assume for the rest of this
section), an anti-isomorphism of algebras with inverse
denoted by $\overline{S}$.  

The reader will note from the axioms above that a Hopf algebra is automatically a weak Hopf algebra.  For a weak Hopf algebra
that is not a Hopf algebra, consider 
a typical groupoid algebra such as $H = M_n(k)$, the $n \times n$-matrices
over $k$ (the groupoid here being a category with $n$ objects
where each $\Hom$-group has a single invertible arrow).
Let $e_{ij}$ denote the $(i,j)$-matrix unit.  
For example, $M_n(k)$ is a weak Hopf algebra
with the counit given by $\eps(e_{ij}) = 1$, comultiplication by $\cop(e_{ij})
= $ $e_{ij} \o e_{ij}$ and antipode given by $S(e_{ij}) = e_{ji}$
for each $i,j = 1,\ldots,n$
(extending the Hopf algebra structure of group algebras).  In this
case, $H^L = H^R$ and is equal to the diagonal matrices.  
The corresponding projections are given by
$\Pi^L(e_{ij}) = e_{ii}$ $= \overline{\Pi}^L(e_{ij})$ and $\Pi^R(e_{ij})= e_{jj}$
$= \overline{\Pi}^R(e_{ij})$. Note that $\eps(1_H) = n 1_k$ which is zero
if the characteristic of $k$ divides $n$.  

There are a number of equations in the subject that we
will need later (cf.\ \cite[2.8, 2.9, 2.24]{BNS}):
\begin{eqnarray}
\Pi^L & = & S \circ \overline{\Pi}^L \\
\Pi^R & = & S \circ \overline{\Pi}^R \label{eq: over2} \\
\overline{S}(a\2)a\1 & = & \overline{\Pi}^R(a) \label{eq: over1} \\
a\2 \overline{S}(a\1) &=& \overline{\Pi}^L(a) \\
a\1 \o \Pi^L(a\2) & = & 1\1 a \o 1\2  \\
\Pi^R(a\1) \o a\2 & = & 1\1 \o a 1\2  \label{eq: pi-are} \\
\Pi^R(a)b & = & b\1 \eps(ab\2)  \\
a\Pi^L(b) & = & \eps(a\1 b) a\2  \label{eq: pi-ell} \ \ \ \ \ (\forall \, a, b \in H) 
\end{eqnarray}
where e.g.\ eq.~(\ref{eq: over1}) follows from applying the inverse-antipode
to eqs.~(\ref{eq: over2}) and~(\ref{eq: bns2}). 

We recall the definition of a right $H$-comodule algebra $A$,
its subalgebra of coinvariants, and Galois coaction 
for $H$ a weak bialgebra (e.g.\ in \cite{CDG}):

\begin{definition}
Let $H$ be a weak bialgebra with $A, H$ both $k$-algebras.
$A$ is a right $H$-comodule algebra if there is a right $H$-comodule
structure $\rho: A \to A \o_k H$ such that $\rho(ab) = \rho(a) \rho(b)$ for each $a,b \in A$ and any of the equivalent
conditions for $\rho(a) := a\0 \o a\1$ are satisfied:
\begin{eqnarray}
1\0 \o 1\1 & \in & A \o H^L \label{eq: cdg2} \\\
a\0 \o \Pi^L(a\1) & = & 1\0 a \o 1\1 \label{eq: cdg1} \\
a\0 \o \overline{\Pi}^R(a\1) & = & a 1\0 \o 1\1 \ \ \ \ (\forall \, a\in H) \label{eq: cdg3} \\
1\0 \o 1\1 \o 1\2 & = & (\rho(1_A) \o 1_H)(1_A \o \cop(1_H)) 
\end{eqnarray}
The coinvariants are defined by 
$$ B := \{ b \in A \| b\0 \o b\1 = 1\0 b \o 1\1  = b 1\0 \o 1\1 \}, $$
the second equation following from equations directly above.
We say $A$ is a weak $H$-Galois extension of $B$ if
the mapping $\beta: A \o_B A \to A \o_k H$ given
by $\beta(a \o a') = a {a'}\0 \o {a'}\1$ is bijective onto
$$ \overline{A \o H} = (A \o H )\rho(1) = \{ a 1\0 \o h 1\1 \|
a \in A, h \in H \}. $$  
\end{definition}

For finite dimensional weak Hopf algebras and their actions,
we only need require $\beta$ be surjective in the definition of weak Hopf-Galois extension,
as $\beta$ is automatically injective by \cite{Bo, BTW} or corollary~\ref{cor-surj} below.  
Note that $\Im \rho \subseteq \overline{A \o H}$,
an $A$-$B$-sub-bimodule
and that $\beta$ is an $A$-$B$-bimodule morphism
w.r.t.\ the structure $a' \cdot (a \o h) \cdot b = a'ab \o h$
on $A \o H$.  
These definitions correspond to the case of a separable base algebra in the definitions of right 
comodule algebras, Galois coring and Galois coactions for bialgebroids
given in \cite{Karl, LK}.
 
We now establish the Hopf algebra analogue of an alternate 
Galois mapping characterizing Galois extension. This would
correspond to working with a left-handed version of
the Galois coring considered in \cite{CDG}. 

\begin{prop}
Suppose $H$ is a weak Hopf algebra and $A$ a right $H$-module
algebra with notation introduced above.  
Let $\beta': A \o_B A \to A \o H$ be defined by 
\begin{equation}
\beta'(a \o a') = a\0 a' \o a\1
\end{equation}
and $\eta: A \o H \to A \o H$ be the map defined by
\begin{equation}
\eta(a \o h) = a\0 \o a\1 S(h).
\end{equation} 
Then $\beta' = \eta \circ \beta$ and $\beta: A \o_B A \to \overline{A \o H}$ is respectively injective, surjective or bijective iff
$\beta'$ is injective, surjective or bijective onto
$$ \overline{\overline{A \o H}} := \rho(1)(A \o H). $$ In particular, 
$A \| B$ is a weak
$H$-Galois extension iff $\beta':A \o_B A \to \overline{\overline{A \o H}}$ is bijective.
\end{prop}
\begin{proof}
Notice that $\overline{\overline{A \o H}}$
is a $B$-$A$-sub-bimodule of $A \o H$,
 and that $\Im \eta$ and $\Im \beta' \subseteq \overline{\overline{A \o H}}$.
Next note that an application of eq.~(\ref{eq: cdg1}) gives
\begin{eqnarray*}
\eta \beta(a \o a') & = & \eta(a{a'}\0 \o {a'}\1) \\
                    & = & a\0 {a'}\0 \o a\1 {a'}\1 S({a'}\2) \\
                    & = & a\0 {a'}\0 \o a\1 \Pi^L({a'}\1) \\
                    & = & a\0 1\0 a' \o a\1 1\1 \\
                    & = & a\0 a' \o a\1 = \beta'(a \o a').
\end{eqnarray*}
We define another linear self-mapping of $A \o H$ given
by $\overline{\eta}(a\o h) = a\0 \o \overline{S}(h)a\1$. 
Note that $\Im \overline{\eta}$ and $\Im \beta \subseteq \overline{A \o H}$.

Let $p: A \o H \to \overline{A \o H}$, $\overline{p}: A \o H \to 
\overline{\overline{A \o H}}$ be the straightforward projections given
by $p(a \o h) = a1\0 \o h1\1$, and $\overline{p}(a \o h) = 1\0 a \o 1\1 h$.
We show below that $\eta \circ p = \eta$, $\overline{\eta} \circ \overline{p} = \overline{\eta}$,
$\eta \circ \overline{\eta} = \overline{p}$ and $\overline{\eta} \circ \eta = p$, from
which it follows that the restrictions of $\eta$, $\overline{\eta}$ to $\overline{A \o H}$, $\overline{\overline{A \o H}}$
are inverses to one another, so that there is a commutative triangle connecting $\beta$, $\beta'$ via
$\eta$.  

$$\begin{diagram}
& & A \o_B A & & \\
&\SW_{\beta'} & & \SE_{\beta} & \\
\overline{\overline{A \o H}} && \lTo^{\cong}_{\eta} && \overline{A \o H} 
\end{diagram}$$

First, we note that $\eta \circ p = \eta$ since
\begin{eqnarray*}
\eta(a 1\0 \o h 1\1) & = & a\0 1\0 \o a\1 1\1 S(h 1\2) \\
                       & = & a\0 1\0 \o a\1 \Pi^L(1\1) S(h) \\
                       & = & a\0 \o a\1 S(h) = \eta(a \o h)
\end{eqnarray*}
by eqs.~(\ref{eq: bns1}) and~(\ref{eq: cdg2}).

Secondly, we note that $\overline{\eta} \circ \overline{p} = \overline{\eta}$ since
\begin{eqnarray*}
\overline{\eta}(1\0 a \o 1\1 h) & = & 1\0 a\0 \o \overline{S}(h) \overline{S}(1\2)1\1 a\1 \\
                                & = & 1\0 a\0 \o \overline{S}(h) \overline{\Pi}^R(1\1) a\1 \\
                                & = & 1\0 a\0 \o \overline{S}(h) 1\1 a\1 = \overline{\eta}(a \o h)
\end{eqnarray*}
by eqs.~(\ref{eq: over1}) and~(\ref{eq: cdg3}).  

Next we note that $\overline{\eta} \circ \eta = p$ since
\begin{eqnarray*}
\overline{\eta}(a\0 \o a\1 S(h)) & = & a\0 \o \overline{S}(a\2 S(h)) a\1 \\
                                 & = & a\0 \o h \overline{\Pi}^R(a\1) \\
                                 & = & a 1\0 \o h 1\1 = p(a \o h)
\end{eqnarray*}
by eqs.~(\ref{eq: over1}) and~(\ref{eq: cdg3}).  

Finally we note $\eta \circ \overline{\eta} = \overline{p}$ since
\begin{eqnarray*}
\eta(a\0 \o \overline{S}(h) a\1) & = &  a\0 \o a\1 S(a\2) h \\
                                & = & a\0 \o \Pi^L(a\1) h = \overline{p}(a\o h)
\end{eqnarray*}
by eq.~(\ref{eq: cdg1}). 
\end{proof}

Again let $H$ be a finite dimensional weak Hopf algebra. Recall that the $k$-dual $H^*$ is also a weak Hopf algebra,
by the self-duality of the axioms, and acts on $H$ by the usual right action $x \leftharpoonup \psi =$ $ \psi(x\1) x\2$
and a similarly defined left action. In addition, a right $H$-comodule algebra $A$ corresponds
to a left $H^*$-module algebra $A$ via $\psi \cdot a := a\0 \psi(a\1)$ \cite{NV}.   
Following Kreimer-Takeuchi and Schneider, there are two proofs that surjectivity of $\beta$ is all
that is needed in the definition of a weak Hopf $H$-Galois extension \cite{Bo, BTW}.  As a corollary of the proposition,
we offer a third and direct proof. 

\begin{cor} \cite{Bo, BTW}
\label{cor-surj}
Let $A$ be a right $H^*$-comodule algebra and $B$ its subalgebra of coinvariants $A^{\rm co \, H^*}$.
If $\beta: A \o_B A \to \overline{A \o H^*}$ is surjective, then the natural module $A_B$ is f.g.\
projective and $A | B$ is a weak $H^*$-Galois extension.
\end{cor}
\begin{proof}
We know from \cite{V} that $H$ and $H^*$ are both Frobenius algebras with nondegenerate left integral 
$t \in H$ satisfying $ht = \Pi^L(h)t$ for all $h \in H$
as well as $t \leftharpoonup T = 1_H$ for some $T \in H^*$.  Since 
$\beta$ is surjective, there are finitely many paired elements $a_i, b_i \in A$ such that
$$ 1\0 \o T 1\1 = \sum_i a_i {b_i}\0 \o {b_i}\1. $$
Let $\phi_i(a) := t \cdot (b_i a)$ for every $a \in A$.  
Then $\{ a_i \}$, $\phi_i$ are dual bases for the module $A_B$ by a computation that $\sum_i a_i \phi_i(a) = a$
for all $a \in A$, 
almost identical with \cite[p.\ 132]{Mo} for Hopf algebra actions (using the identity
$1\0 a\0 \o 1\1 a\1 = a\0 \o a\1$ at one point).

Finally, one shows that $\beta'$ is injective, for if $\sum_j u_j \o v_j \in \ker \beta'$, we
compute $$\sum_j u_j \o v_j = \sum_{i,j} a_i \o \phi_i(u_j)v_j = \sum_{i,j} a_i \o (t\1 \cdot b_i) {u_j}\0 v_j \bra {u_j}\1 , t\2 \ket = 0 $$
as in \cite[p.\ 132]{Mo}. By the proposition, $\beta$ is then
injective, whence a bijection of $A \o_B A$ onto $\overline{A \o H}$. 
\end{proof}

We next offer an example of weak Hopf-Galois extension with an
alternative proof.  For example, if $H = M_n(k)$ considered above,
the Galois map $\beta = $ \newline
$(\mu \o \id)\circ (\id \o \cop)$ given by  $\beta(e_{ij} \o e_{jk}) =$ $e_{ik} \o e_{jk}$
with coinvariants $H^L$ the diagonal matrices and $1\1 \o 1\2 =$
$\sum_i e_{ii} \o e_{ii}$, is an isomorphism by a dimension count. The
general picture is the following: 

\begin{cor} \cite[2.7]{CDG}
\label{cor - olary}
Define a coaction on $H$ by $a\0 \o a\1 = a\1 \o a\2$
for $a \in H$.
Then $H$ is a weak Hopf $H$-Galois extension of its coinvariants
$H^L$.
\end{cor}
\begin{proof}
We note that $H^L \subseteq H^{\rm co \, H}$ since
$\cop(x^L) = 1\1 x^L \o 1\2$ for $x^L \in H^L$ \cite[2.7a]{BNS}.  The converse follows
from $x \in H^{\rm co \, H}$ implies $$ x = \eps(x\1) x\2
= \eps(1\0 x) 1\1 \in H^L.$$

Next we note that $\beta'$ factors into isomorphisms in the following
commutative diagram, where $\sigma: H \o H \to H \o H$ is
the standard twist involution:

$$\begin{diagram}
H \o_{H^L} H &&  \rTo^{\beta'}&& \overline{\overline{H \o H}} \\
\dTo_{\cong}^{q} && && \uTo^{\cong}_{\sigma \circ (S \o S)}    \\
\overline{\overline{H \o H}} &&  \rTo_{\overline{\eta}}^{\cong}   && \overline{H \o H}
\end{diagram}$$

where $q(x \o y) :=$ $ \overline{p}(\overline{S}(x) \o_k y)$
is well-defined since $S(1\1) \o 1\2$ is a separability element
for the separable $k$-algebra $H^L$ \cite[prop.\ 2.11]{BNS}.  Its inverse
is given by $q^{-1}(\overline{p}(x \o y))$ $= S(x) \o y$. 
The mapping $\sigma \circ (S \o S)$ has an obvious inverse
and is well-defined since $S(1) = 1$ and $S$ is an anti-coalgebra
homomorphism.  
\end{proof}

We provide the complete proof that a weak Hopf $H$-Galois extension
is depth two \cite[3.2]{LK}:
\begin{cor}
A weak $H$-Galois extension $A \| B$ is right and left depth two.
\end{cor}
\begin{proof}
The algebra extension $A \| B$ is right D2 since the Galois mapping $\beta: A \o_B A \stackrel{\cong}{\longrightarrow}$
$\overline{A \o H}$ and the projection $p: A \o H \to \overline{A \o H}$ are $A$-$B$-bimodule morphisms \cite[3.1]{Karl}.
Whence $A \o_B A$ is $A$-$B$-isomorphic to a direct summand $\Im p$ within $\oplus^n A$ where $n = \dim H$.

Similarly $A \| B$ is left D2 since the alternate Galois isomorphism $\beta'$  and projection $\overline{p}$
onto $\overline{\overline{A \o H}}$ 
are both $B$-$A$-bimodule morphisms.
\end{proof}

The proof of the corollary sidesteps the problem of showing $A \| B$ is a Frobenius extension,
which then implies that left D2 $\Leftrightarrow$ right D2.  It is likely that a
direct proof using $\beta$ and $\beta'$ that a weak $H$-Galois extension is Frobenius may be given 
since there are nondegenerate integrals in $H^*$ which would define a Frobenius homomorphism via the dual
action of $H^*$ on $A$ (with invariants $B$). 
In addition we have avoided starting only with a weak bialgebra
having Galois action on $A$ and showing  the existence of an antipode on $H$ 
 in extension of \cite{Sch} for Galois actions of bialgebras. 
If we denote 
\begin{equation}
\label{eq: notation}
\sum_i \ell_i(h) \o_B r_i(h) := \beta^{-1}(1\0 \o h 1\1),
\end{equation}
we note that
\begin{eqnarray}
1\0 \o 1\1 S(h) & = & \eta(1 \o h) \\
                & = & \beta' (\beta^{-1}(1\0 \o h 1\1)) \\
                & = & \sum_i {\ell_i(h)}\0 r_i(h) \o {\ell_i(h)}\1,
\end{eqnarray}
which can conceivably be made to descend to a formula for the antipode of $H$ in terms
of just the isomorphism $\beta$.  

We  then propose  two problems and provide some evidence for each.

\begin{prob}
If $H$ is a finite dimensional weak Hopf algebra and $A \| B$ is $H$-Galois, provide a direct proof that 
$A \| B$  is a Frobenius extension (cf.\ \cite[3.7]{BaS}).
\end{prob}

For example, if $H$ is a Galois extension of $H^L$ as in corollary~\ref{cor - olary}
we expect such a Frobenius extension based on Pareigis's theorem that a Frobenius subalgebra $B$ of a Frobenius algebra $A$,
where the natural module $A_B$ is f.g.\ projective and the Nakayama automorphism of $A$ stabilizes $B$,
yields a $\beta$-Frobenius extension $A \| B$ where $\beta$ is the relative Nakayama automorphism of $A$ and $B$ (restrict
one and compose with the inverse of the other).  For example, if $H$ has an $S$-invariant nondegenerate integral,
the Nakayama automorphism is $S^2$ \cite[3.20]{BNS}, also the Nakayama automorphism of $H^L$, so $\beta = \id$
and $H \| H^L$ is an ordinary Frobenius extension.  

\begin{prob}
\label{conj: b} 
If $H$ is a weak bialgebra and $A \| B$ is $H$-Galois, is 
$H$ necessarily a weak Hopf algebra? 
\end{prob}

Again this is true in the special case of the weak Hopf-Galois extension in corollary~\ref{cor - olary},
a result in \cite[Brzezinski-Wisbauer]{BW}; we give another proof which may extend to  the general problem.
Note that the definition of weak Hopf-Galois extension does not make use of an antipode
nor does $H^{\rm co \, H} = H^L$ in corollary~\ref{cor - olary}.

\begin{theorem} \cite[36.14]{BW}
Let $H$ be a weak bialgebra. If the right $H$-coalgebra
$H$ with coaction $\varrho = \cop_H$ is Galois over $H^L$, then
$H$ is a weak Hopf algebra.
\end{theorem}

\begin{proof}
In terms of the notation in eq.~(\ref{eq: notation}) we define an antipode $S: H \to H$ by
\begin{equation}
S(h) = \sum_i \eps({\ell_i(h)}\1 r_i(h)) {\ell_i(h)}\2
\end{equation}

Note that by eq.~(\ref{eq: pi-ell}), $S(h) =$ $ \sum_i \ell_i(h) \Pi^L(r_i(h))$ for $h\in H$.
In order to prove that $S$ satisfies the three eqs.~(\ref{eq: bns2}),~(\ref{eq: bns1}) and~(\ref{eq: bnsa}),
we  note the three equations below for a general right $H$-comodule algebra $A$ over a weak bialgebra $H$
where $A$ is $H$-Galois over its coinvariants $B$;
the proofs are quite similar to those in \cite{Sch}.
\begin{eqnarray}
\sum_i \ell_i(h) \o {r_i(h)}\0 \o {r_i(h)}\1 & = & \sum_i \ell_i(h\1) \o r_i(h\1) \o h\2 \label{eq: s1}  \\
\sum_i a\0 \ell_i(a\1) \o_B r_i(a\1) & = & 1 \o_B a  \label{eq: s2} \\
\sum_i \ell_i(h) r_i(h) & =  & 1\0 \eps(h 1\1 ) \ \ \ \ (\forall \, a \in A, h \in H) \label{eq: s3}
\end{eqnarray}

Next we note three equations in $A \o H$, two of which we need here (and  all three might play a role in an answer to problem~\ref{conj: b}). 
\begin{eqnarray}
\sum_i {\ell_i(h\1)}\0 r_i(h\1 ) \o {\ell_i (h\1 )}\1 h\2 & = & 1\0 \o 1\1 \Pi^R(h) \label{eq: ex1} \\
\sum_i {\ell_i(h\2)}\0 r_i(h\2) \o h\1 {\ell_i(h\2 )}\1 & = & 1\0 \o \Pi^L(h 1\1 ) \label{eq: ex2}
\end{eqnarray}

\begin{equation}
\label{eq: ex3}
\sum_i {\ell_i(h\1 )}\0 r_i(h\1 ) {\ell_i(h\3)}\0 r_i (h\3 ) \o {\ell_i(h\1 )}\1 h\2 {\ell_i(h\3 )}\1 = 
\end{equation} 
$$ \sum_i {\ell_i(h)}\0 r_i(h) \o {\ell_i(h)}\1. $$ 
They are established somewhat similarly to \cite{Sch} and left as exercises.

Applying eq.~(\ref{eq: ex1}) with $A = H$ and $a\0 \o a\1 = a\1 \o a\2$, we obtain one of the antipode axioms:
\begin{eqnarray*}
S(h\1) h\2 & = & \sum_i \eps({\ell_i(h\1)}\1 r_i(h\1 )){\ell_i (h\1 )}\2 h\2  \\
           & = & \eps(1\1 ) 1\2 \Pi^R(h) = \Pi^R(h). \ \ \ \ (\forall \, h \in H) 
\end{eqnarray*}
 
Applying eq.~(\ref{eq: ex2}), we obtain
\begin{eqnarray*}
h\1 S(h\2) & = & \sum_i \eps({\ell_i(h\2)}\1 r_i(h\2) )h\1 {\ell_i(h\2 )}\2 \\
           & = & \eps(1\1 ) \Pi^L(h 1\2 ) = \Pi^L(h) \ \ \ \ (\forall \, h \in H) 
\end{eqnarray*}

Finally we see $S$ is an antipode from the just established eq.~(\ref{eq: bns2}) and applying eq.~(\ref{eq: pi-are}): 
\begin{eqnarray*}
\Pi^R(h\1) S(h\2) & = & \sum_i \Pi^R(h\1) \ell_i(h\2) \Pi^L(r_i(h\2)) \\
                   &=&  \sum_i 1\1 \ell_i(h1\2) \Pi^L(r_i(h1\2 )) \\
                   & = & \sum_i \ell_i(h) \Pi^L(r_i(h)) = S(h)  
\end{eqnarray*}
where we use the general fact that $\beta$ is left $A$-linear, so $\sum_i 1\0 \ell_i(h 1\1) \o r_i(h 1\1) = $
$\beta^{-1}({1'}\0 1\0 \o h {1'}\1 1\1) = $ $ \sum_i \ell_i(h) \o r_i(h)$.  
\end{proof}

%If antipodes
%exist for weak bialgebras with Galois coaction, certain depth two extensions 
%are automatically Frobenius extension.  For example,  
 
%%%%%%%%%%%%%%%%%%%%%%%%%%%%%%%%%%%%%%%%%%%%%%%%%%%%%%%%%%%%%%%%%%%%%%
%%%%%%%%%%%%%%%%%%%%%%%%%%%%%%%%%%%%%%%%%%%%%%%%%%%%%%%%%%%%%%%%%%%%%%
\section{Appendix}
In this section we answer some natural questions about the theory of one-sided depth two extensions.
One of the apparent questions after a reading of proposition~\ref{prop-tee} would
be if the  endomorphism algebra $S$ is also a bialgebroid over the centralizer, to which the next proposition 
provides an answer in the affirmative.
\begin{prop}
\label{prop-ess}
Suppose $A \| B$ is either a right or a left D2 extension with centralizer $R$.
Then $S$ is a left $R$-bialgebroid, which is either
right f.g.\ $R$-projective or left f.g.\ $R$-projective
respectively.
\end{prop}
\begin{proof}
The algebra structure comes from the usual composition of endomorphisms in $S =
\End {}_BA_B$.  The source and target mappings are $s_L(r) = \lambda(r)$
and $t_L(r) = \rho(r)$, whence the structure ${}_RS_R$ is given
by $$ r \cdot \alpha \cdot r' = \lambda(r) \rho(r') \alpha = r \alpha(-) r'.
$$

Suppose now we are given a right D2 structure on $A \| B$ by quasibases
$u_j \in T$, $\gamma_j \in S$.  
The $R$-coring structure on ${}_RS_R$ is given by a coproduct $\cop: S  \to S\o_R S$ defined by 
\begin{equation}
\label{eq: coproduct}
\cop(\alpha) = \sum_j \gamma_j \o u^1_j \alpha(u^2_j -) ,
\end{equation}
and a counit $\eps: S \to R$ given by
\begin{equation}
\label{eq: counit}
\eps(\alpha) = \alpha(1_A)
\end{equation}
Clearly $\eps$ is an $R$-$R$-bimodule mapping with $\eps(1_S) = 1_A$,
satisfying the counitality equations and $$\eps(\alpha \beta) = \eps(\alpha s_L(\eps(\beta))) = \eps(\alpha t_L(\eps(\beta))).$$ Also $\cop$ is right $R$-linear and $\cop(1_S) = 1_S \o_R 1_S$. 
By
making the identification 
$$ S \o_R S \cong \Hom ({}_B A \o_B A_B, {}_B A_B), \ \ \ \alpha \o \beta \longmapsto
(a \o a' \mapsto \alpha(a) \beta(a') )$$
with inverse $F \mapsto \sum_j \gamma_j \o u^1_j F(u^2_j \o - )$, 
we see that the coproduct is left $R$-linear, satisfies $\alpha\1 t_L(r) \o \alpha\2 =$ $ \alpha\1 \o \alpha\2 s_L(r)$
for all $r \in R$, and $\cop(\alpha \beta) = \cop(\alpha)\cop(\beta)$ for all $\alpha, \beta \in S$.  For with the independent variables $x,x' \in A$, $\alpha, \beta \in S$
and $r \in R$,
each of these expressions becomes equal  to $ r\alpha(xx')$, $\alpha(xrx')$, and $ \alpha(\beta(xx'))$
respectively. 

The coproduct $\cop$ is coassociative since 
$$ S \o_R S \o_R S \stackrel{\cong}{\longrightarrow} \Hom ({}_BA \o_B A \o_B A_B, {}_BA_B), \ 
\alpha \o \beta \o \gamma \longmapsto ( x \o y \o z \mapsto \alpha(x)\beta(y) \gamma(z) ) $$
with inverse given by
\begin{equation}
F \mapsto \sum_{i,j,k} \gamma_i \o u^1_i \gamma_j(u^2_i \gamma_k(-)) \o u^1_j F(u^2_j u^1_k \o u^2_k \o -) 
\end{equation}
Applying this identification to $(\cop \o \id_S) \cop(\alpha)$ and to $(\id_S \o \cop)\cop(\alpha)$
then to $x \o_B y \o_B z$ both expressions equal $\alpha(xyz)$. 

$S_R$ is f.g.\ projective since for each $\alpha \in S$, we have $\alpha = \sum_j \gamma_j h_j(\alpha)$
where $h_j \in \Hom (S_R, R_R)$ is defined by $h_j(\alpha) = u_j^1 \alpha(u_j^2)$. 

The proof that given left D2 quasibases $t_i \in T$, $\beta_i \in S$,
we have left f.g.\ projective left bialgebroid $S$ with identical bialgebroid structure is similar
and therefore omitted.
\end{proof}

Suppose $A \| B$ is right D2.  Then we have seen that $S$ is a right finite projective left bialgebroid
while $T$ is a left finite projective right bialgebroid.  There is a nondegenerate pairing between
$S$ and $T$ with values in the centralizer $R$ given by $\bra t \| \alpha \ket := $ $t^1 \alpha(t^2)$,
since \begin{equation}
\label{eq: pair}
\eta: {}_RT \stackrel{\cong}{\longrightarrow} \Hom (S_R, R_R) 
\end{equation}
via $\eta(t) = \bra t \| - \ket$ with inverse $\phi \mapsto$ $\sum_j \phi(\gamma_j)u_j$. 
By proposition \cite[2.5]{KS} a right f.g.\ projective left bialgebroid $S$ has a right $R$-bialgebroid
$R$-dual $S^*$.  The question is then if the  bialgebroid $S^*$ is isomorphic to the bialgebroid $T$
via $\eta$?  The question is partly answered in the affirmative by corollary \cite[5.3]{KS},
where it is shown without using left D2 quasibases that $T$ and $S^*$ are isomorphic via the pairing above
as algebras and $R$-$R$-bimodules.
\begin{cor}
Suppose $A \| B$ is right D2.  Then $T$ is isomorphic as right bialgebroids over $R$ to the right
$R$-dual of $S$ via $\eta$.  If $A \| B$ is left D2, then $T$ is isomorphic to
the bialgebroid left $R$-dual of $S$.
\end{cor}
\begin{proof}
What remains to check in the first statement is that $\eta$ is a homomorphism of $R$-corings
using right D2 quasibases.
We compute:
\begin{eqnarray*} 
\bra t\1 \cdot \bra t\2 \| \alpha' \ket \| \alpha \ket & = & \sum_j \bra t^1 \o_B \gamma_j(t^2)u^1_j \alpha'(u^2_j) \| \alpha \ket \\
& = & t^1 \alpha (\alpha'(t^2)) = \bra t \| \alpha \circ \alpha' \ket,
\end{eqnarray*}
Whence $\cop(\eta(t)) = \eta(t\1) \o \eta(t\2)$ by uniqueness \cite[2.5 (41)]{KS}. 

The proof of the last statement is similar to the first in using the pairing $[\alpha \| t ] :=$ $ \alpha(t^1)t^2$
and the right bialgebroid of the left dual of a left bialgebroid in \cite[2.6]{KS}. The details are left
to the reader.
\end{proof}

%%%%%%%%%%%%%%%%%%%%%%%%%%%%%%%%%%%%%%%%%%%%%%%%%%%%%%%%%%%%%%%%%%%%%

\end{document}